\documentclass[reqno,12pt]{amsart}
\usepackage[T1]{fontenc}
\usepackage{lmodern}
\usepackage{amsmath,amssymb}
\usepackage{microtype}
\usepackage{needspace}
\usepackage[colorlinks=true,linkcolor=blue,citecolor=blue,urlcolor=blue]{hyperref}
\hypersetup{pdftitle={Uniformly positive entropy and mixing in minimal systems},
	pdfauthor={Shuhao Zhang and Xiaomin Zhou}}
\makeatletter
\@namedef{subjclassname@2020}{\textup{2020} Mathematics Subject Classification}
\makeatother
\numberwithin{equation}{section}
\newtheorem{thm}{Theorem}[section]
\newtheorem{prop}[thm]{Proposition}
\newtheorem{lem}[thm]{Lemma}
\newtheorem{cor}[thm]{Corollary}
\newtheorem{conj}[thm]{Conjecture}
\theoremstyle{definition}
\newtheorem{problem}[thm]{Problem}
\theoremstyle{remark}
\newtheorem{rem}[thm]{Remark}
\newcommand{\N}{\mathbb N}
\newcommand{\Z}{\mathbb Z}
\newcommand{\Asy}{\operatorname{Asy}}
\newcommand{\diam}{\operatorname{diam}}
\newcommand{\dist}{\operatorname{dist}}

\begin{document}
	\title[Uniformly positive entropy and mixing]{Uniformly positive entropy and mixing in minimal systems}
	\author[S. Zhang]{Shuhao Zhang}
	\address[S. Zhang]{School of Mathematical Sciences, University of Science and Technology of China, Hefei, Anhui 230026, P.R. China}
	\email[S. Zhang]{yichen12@mail.ustc.edu.cn}
	\author[X. Zhou]{Xiaomin Zhou}
	\address[X. Zhou]{School of Mathematics and Statistics, Huazhong University of Science and Technology, Wuhan, Hubei 430074, China}
	\email[X. Zhou]{zxm12@mail.ustc.edu.cn}
	\date{}
	
	\begin{abstract}
		Every minimal system with uniformly positive entropy of order two is
		topologically mixing. In a minimal system, mixing along an infinite
		set $A\subset\N_0$ is equivalent to $\dist(T^nU,T^nV)\to0$ along $A$
		for all nonempty open $U,V$, and to unbounded complexity of each
		standard binary open cover along every infinite subset of $A$.
		For $A=\N_0$ it is also equivalent to $\Delta^*$-transitivity.
		A local distance condition characterizes finite cyclic decompositions
		with mixing return maps. Such a decomposition exists precisely when,
		for each infinite $S\subset\N_0$, the closure of some proximal cell
		along $S$ has positive measure for some invariant probability.
		In a minimal weakly mixing system, for fixed infinite $S$, these cells
		are either all dense or all have closures of measure zero for every
		invariant probability. We also give proximal-cell criteria for
		Furstenberg families and mild mixing.
	\end{abstract}
	\maketitle
	
	\section{Introduction}
	
	Blanchard \cite{B92} introduced uniformly positive entropy (UPE),
	which requires positive entropy for every open cover by two non-dense sets.
	A system has UPE if and only if every pair of distinct points is an
	entropy pair \cite{B93}.

	For invertible probability-preserving systems, positive entropy for every
	nontrivial two-atom partition implies the Kolmogorov property and hence
	strong mixing \cite{GY09}.
	Blanchard showed that UPE implies weak mixing and gave nonminimal
	examples in which topological mixing fails \cite{B92,B93}. Huang and Ye
	showed that UPE implies mild mixing \cite{HY06}.

	Glasner and Weiss \cite[Theorem~C]{GW94} proved that every ergodic positive-entropy process has a strictly ergodic UPE model. Huang, Shao and Ye \cite[Theorem~5.10]{HSY05} proved that a minimal system with UPE of all finite orders, equivalently a minimal topological $K$-system, is mixing. For general systems, Huang and Ye \cite{HY06} showed that UPE of order $n$ need not imply UPE of order $n+1$.
	
	Glasner and Weiss \cite[Section~5, Remark~2]{GW94} asked whether every
	minimal system with UPE is topologically mixing. Huang, Shao and Ye
	restated the question in \cite[Problem~1]{HSY05}. It also appears in
	\cite[Section~16]{GY09}.

	\begin{thm}\label{thm:upe}
		Every minimal system with uniformly positive entropy is topologically mixing.
	\end{thm}

	Throughout, a dynamical system $(X,T)$ consists of a nonempty compact metric space $(X,\rho)$ and a continuous map $T:X\to X$. Minimality means that every forward orbit is dense; in particular, a minimal map is onto. Invertibility is assumed only where it is stated explicitly.
	Unless otherwise specified, UPE refers to UPE of order two and mixing to
	topological mixing.
	
	To prove Theorem~\ref{thm:upe}, we consider the forward asymptotic relation
	\[
	\Asy(T)=\{(x,y)\in X^2:\rho(T^nx,T^ny)\to0\text{ as }n\to\infty\}.
	\]
	For homeomorphisms, Blanchard, Host and Ruette
	\cite[Proposition~4]{BHR02} proved that every entropy pair belongs to
	$\overline{\Asy(T)}$. Thus a UPE homeomorphism satisfies
	$\overline{\Asy(T)}=X^2$. For nonempty subsets $E,F\subset X$, write
	$\dist(E,F)=\inf\{\rho(x,y):x\in E,\ y\in F\}$.
	Given nonempty open sets $U,V\subset X$, choose
	$(x,y)\in\Asy(T)\cap(U\times V)$. Then
	\[
	\dist(T^nU,T^nV)\leq\rho(T^nx,T^ny)\longrightarrow0.
	\]
	We show that for a minimal homeomorphism this distance condition implies mixing.
	Indeed, if mixing fails, there are nonempty open sets $U_0,V_0$ and times
	$n_j\to\infty$ such that $T^{n_j}U_0\cap V_0=\varnothing$ for every $j$.
	After passing to a subsequence as in Section~\ref{sec:criterion}, let
	$L(x)$ be the set of all limits of $T^{n_{j_k}}x_k$ with $x_k\to x$.
	The distance condition gives $L(x)\cap L(y)\ne\varnothing$ for all $x,y\in X$.
	Proposition~\ref{prop:alllimits} first deduces weak mixing from these
	pairwise intersections. In the induction from $s$ points to $s+1$,
	Lemma~\ref{lem:completion} and minimality give a common point of
	$L(x),L(T^dx),\ldots,L(T^{sd}x)$ for every $x\in X$ and $d\ge1$.
	By minimality and weak mixing, the tuples $(x,T^dx,\ldots,T^{sd}x)$ are
	dense in $X^{s+1}$. Since the set of tuples with this intersection
	property is closed, the induction gives the finite intersection property
	for $\{L(x):x\in X\}$.
	Its intersection is nonempty by compactness and is $T$-invariant since
	$L(Tx)=TL(x)$. By minimality it equals $X$, contradicting
	$L(x)\cap V_0=\varnothing$ for $x\in U_0$.
	For a noninvertible minimal system with UPE, we apply the homeomorphism
	case to its natural extension, which is again minimal and has UPE.
	Mixing passes to the original system (Lemma~\ref{lem:naturalextension}).
	
	We write $\N=\{1,2,\ldots\}$ and $\N_0=\N\cup\{0\}$.
	For an infinite set $A\subset\N_0$, we say that $(X,T)$
	is \emph{mixing along $A$} if, for every pair of nonempty open sets
	$U,V\subset X$, we have $T^nU\cap V\ne\varnothing$ for all sufficiently
	large $n\in A$.
	
	Blanchard, Host and Maass \cite{BHM00} studied open-cover complexity
	and its relation to scattering. Huang and Ye \cite{HY04} showed that
	standard binary open covers suffice to characterize scattering. They
	also proved that for minimal systems full-scattering is equivalent to
	mixing. Full-scattering requires unbounded complexity along every
	infinite sequence for each finite open cover whose members are non-dense.
	Here we characterize mixing along a prescribed infinite set $A$ by the
	complexity of standard binary open covers along the infinite subsets of $A$.
	
	For a finite open cover $\alpha$ and a finite nonempty set
	$F\subset\N_0$, write
	\[
	c_\alpha(F)=N\left(\bigvee_{n\in F}T^{-n}\alpha\right),
	\]
	where $N(\mathcal U)$ is the least cardinality of a subcover.
	A binary open cover is \emph{standard} if neither member is dense.
	
	For an infinite set $B\subset\N_0$, let
	$\Delta(B)=\{c-b:b,c\in B,\ b<c\}$. A system is
	\emph{$\Delta^*$-transitive} if every open-set return-time set meets
	every such $\Delta(B)$, and \emph{$\Delta^*$-mixing} if its Cartesian
	square is $\Delta^*$-transitive.
	
	\begin{thm}\label{thm:criterion}
		Let $(X,T)$ be minimal and let $A\subset\N_0$ be infinite.
		The following conditions are equivalent:
		\begin{enumerate}
			\item $(X,T)$ is mixing along $A$.
			\item For all nonempty open sets $U,V\subset X$,
			\begin{equation}\label{eq:distance}
				\dist(T^nU,T^nV)\longrightarrow0
				\qquad(n\to\infty,\ n\in A).
			\end{equation}
			\item For every standard binary open cover $\alpha$ and every
			infinite $S=\{s_1<s_2<\cdots\}\subset A$,
			\[
			c_\alpha(\{s_1,\ldots,s_m\})\longrightarrow\infty.
			\]
		\end{enumerate}
		When $A=\N_0$, these conditions are also equivalent to:
		\begin{enumerate}
			\setcounter{enumi}{3}
			\item $(X,T)$ is $\Delta^*$-transitive.
			\item $(X,T)$ is $\Delta^*$-mixing.
		\end{enumerate}
		If $(X,T)$ is not mixing along $A$, there are a standard binary
		open cover $\alpha$ and an infinite $B\subset A$ such that
		\begin{equation}\label{eq:introbinary}
			c_\alpha(F)=2
			\qquad\text{for every nonempty finite }F\subset B.
		\end{equation}
	\end{thm}
	For a system with UPE, the case $n=2$ of \cite[Theorem~7.1]{HLY12}
	gives positive sequence entropy for every standard binary open cover
	along every infinite time set. This implies~(3). UPE also implies
	$\Delta^*$-transitivity by \cite[Corollary~7.3]{HLY12}. Thus either
	$(3)\Rightarrow(1)$ or $(4)\Rightarrow(1)$ gives Theorem~\ref{thm:upe}
	when $A=\N_0$.

	\begin{rem}\label{rem:mixingconsequences}
		Condition~(3) of Theorem~\ref{thm:criterion} holds if and only if
		for every standard binary open cover $\alpha$ and every infinite
		$S\subset A$, there is an infinite $B\subset S$ such that
		\[
		c_\alpha(F)=2^{|F|}
		\qquad(\varnothing\ne F\subset B\text{ finite}).
		\]
	\end{rem}

	Huang and Ye \cite[Question~1]{HY04} asked whether $\Delta^*$-mixing
	implies full-scattering; see also \cite[Section~3.1, p.~9]{DA17}.
	Xu \cite[Theorem~1.2]{Xu26} gave a negative answer in the nonminimal
	setting. For minimal systems, however, Theorem~\ref{thm:criterion}
	shows that even $\Delta^*$-transitivity implies mixing, and hence
	full-scattering by \cite[Theorem~5.8]{HY04}.

	The distance condition in Theorem~\ref{thm:criterion}
	is global. We now consider its local version.
	A nontrivial finite cyclic factor obstructs mixing even when
	the return map on each fiber is mixing. When the
	maximal equicontinuous factor of a minimal system
	is finite, the return map on each fiber is weakly mixing
	\cite[Proposition~29 and the proof of Proposition~30]{BHS08}.
	
	Following Huang, Shao and Ye \cite{HSY04}, for an infinite set
	$S\subset\N_0$ and $x\in X$ we define the \emph{proximal cell of $x$ along $S$} by
	\begin{equation}\label{eq:introcell}
		P_S(x)=\left\{y\in X:
		\liminf_{\substack{n\to\infty\\n\in S}}
		\rho(T^nx,T^ny)=0\right\}.
	\end{equation}
	In particular, $P_{\N_0}(x)$ is the usual proximal
	cell. Furstenberg \cite{F81} proved that these cells
	are dense $G_\delta$ sets in minimal weakly mixing
	systems, and Akin and Kolyada \cite[Theorem~3.8]{AK03}
	extended this result to weakly mixing systems.
	Huang, Shao and Ye \cite[Theorem~4.2(1)]{HSY04}
	proved that, for a system admitting an invariant
	measure of full support, mixing is equivalent to
	the density of $P_S(x)$ for every $x\in X$ and every
	infinite $S\subset\N_0$.
	
	Let $M_T(X)$ denote the set of $T$-invariant Borel probability measures on $X$.
	For a minimal system, consider the weaker requirement that for each infinite
	$S\subset\N_0$ there exist $x_S\in X$ and $\mu_S\in M_T(X)$ with
	$\mu_S(\overline{P_S(x_S)})>0$. By the following theorem, this is
	equivalent both to a local distance condition and to the existence of a
	finite cyclic decomposition with mixing return maps.
	
	\begin{thm}\label{thm:local}
		For a minimal system $(X,T)$, the following conditions are equivalent:
		\begin{enumerate}
			\item There is a nonempty open set $O\subset X$ such that
			\begin{equation}\label{eq:localdistance}
				\dist(T^nU,T^nV)\longrightarrow0
				\qquad(n\to\infty)
			\end{equation}
			for all nonempty open sets $U,V\subset O$.
			\item There are $x\in X$ and a nonempty open set $O\subset X$
			such that \eqref{eq:localdistance} holds for every open
			neighborhood $U$ of $x$ and every nonempty open set $V\subset O$.
			\item There is a finite decomposition into nonempty clopen sets
			\begin{equation}\label{eq:cyclic}
				X=\bigsqcup_{i\in\Z/p\Z}C_i,\qquad TC_i=C_{i+1},
			\end{equation}
			such that each $T^p|_{C_i}$ is minimal and mixing.
			\item For every infinite $S\subset\N_0$, there is $x_S\in X$ such that
			\[
			\operatorname{int}\overline{P_S(x_S)}\ne\varnothing.
			\]
			\item For every infinite $S\subset\N_0$, there are $x_S\in X$
			and $\mu_S\in M_T(X)$ such that
			\begin{equation}\label{eq:intropositivecell}
				\mu_S\bigl(\overline{P_S(x_S)}\bigr)>0.
			\end{equation}
		\end{enumerate}
		The sets $C_i$ are the fibers of the maximal equicontinuous factor
		and are unique up to cyclic relabeling. If $T$ is weakly mixing
		or totally minimal, or if $X$ is connected, then $p=1$.
	\end{thm}
	
	\begin{rem}\label{rem:localvariants}
		\begin{enumerate}
			\item If $\overline{\Asy(T)}$ or the closure of some stable class has
			nonempty interior, then the cyclic decomposition in (3) exists
			(Corollary~\ref{cor:localasy}). Under this decomposition, return-time
			sets are eventually determined by residue classes modulo $p$
			(Corollary~\ref{cor:residues}).
			
		\item Conditions~(3)--(5) remain equivalent when ``mixing'' is replaced by ``mildly mixing'' in (3) and ``infinite set'' by ``IP-set'' in (4) and (5).
			
		\item Under either decomposition, we have
		\[
		\overline{P_S(x)}=C_i\qquad(x\in C_i)
		\]
		for every infinite $S$ in the mixing case and every IP-set $S$ in the mildly mixing case. The conclusions of
		Theorem~\ref{thm:local} about the maximal equicontinuous factor and the
		case $p=1$ also hold for the mildly mixing decomposition.
		
		Conversely, if the corresponding decomposition does not exist, there is an infinite set $S$, respectively an IP-set $S$, such that
		\[
		\mu(\overline{P_S(x)})=0
		\qquad(x\in X,\ \mu\in M_T(X)).
		\]
		In conditions~(4) and (5), the choices of $x_S$ and $\mu_S$ may depend on $S$.
		\end{enumerate}
	\end{rem}
	
	Conditions~(4) and (5) of Theorem~\ref{thm:local} concern proximal
	cells along every infinite time set. In a minimal weakly mixing system,
	fix an infinite set $S\subset\N_0$. If some invariant probability measure gives
	positive mass to the closure of one proximal cell along $S$, then
	$P_S(x)$ is dense for every $x\in X$.
	\begin{thm}\label{cor:celldichotomy}
		Let $(X,T)$ be minimal and weakly mixing, and let
		$S\subset\N_0$ be infinite. The following conditions are equivalent:
		\begin{enumerate}
			\item $P_S(x)$ is a dense $G_\delta$ subset of $X$ for every $x\in X$.
			\item There is $x\in X$ such that
			$\operatorname{int}\overline{P_S(x)}\ne\varnothing$.
			\item There are $x\in X$ and $\mu\in M_T(X)$ such that
			$\mu(\overline{P_S(x)})>0$.
		\end{enumerate}
		If these conditions fail, then
		\[
		\mu(\overline{P_S(x)})=0
		\qquad(x\in X,\ \mu\in M_T(X)).
		\]
	\end{thm}
	
	\begin{rem}\label{rem:proximalextensions}
		\begin{enumerate}
			\item In the nondense case, visits to cell closures satisfy the uniform estimate
			\[
			\lim_{N\to\infty}\sup_{x,y\in X}
			\frac1N\bigl|\{0\le j<N:
			T^jy\in\overline{P_S(x)}\}\bigr|=0.
			\]
			
			\item If $(X,T)$ is mixing, the dense case holds for every infinite
			$S\subset\N_0$. If a minimal weakly mixing system fails to be mixing,
			the null case holds for some infinite $S$. For $S=\N_0$, weak mixing
			suffices for the dense case.

		\end{enumerate}
	\end{rem}
	For a minimal system $(X,T)$ and a full Furstenberg family $\mathcal F$,
	the criterion of Huang, Shao and Ye \cite[Proposition~3.3]{HSY04} for
	$\mathcal F$-mixing requires regional proximality of all tuples of every
	order $n\ge2$ along each set in the dual family.
	Theorem~\ref{thm:family} shows that pairs suffice for every proper
	Furstenberg family $\mathcal F$ containing all cofinite sets.
	If $(X,T)$ is also weakly mixing, it suffices that for each $S$ in the dual family
	there is a point $x$ such that $\overline{P_S(x)}$ has nonempty
	interior, or equivalently that for each such $S$ some cell closure has
	positive measure for an invariant probability measure.

	After the preliminaries in Section~\ref{sec:preliminaries}, we prove the
	distance criterion and deduce Theorem~\ref{thm:upe} in Section~\ref{sec:criterion}.
	Section~\ref{sec:binary} uses binary open covers to complete the proof of
	Theorem~\ref{thm:criterion}. Section~\ref{sec:local} derives the cyclic
	decomposition from the local distance condition. Section~\ref{sec:families}
	proves the proximal-cell criteria and the mildly mixing variant in
	Remark~\ref{rem:localvariants}. In Section~\ref{sec:questions}, a Sturmian
	extension gives a negative answer to the questions of Kamiński, Siemaszko
	and Szymański on relative weak mixing. We also discuss asymptotic pairs
	of measures and higher-order mixing.
	
	\section{Preliminaries}\label{sec:preliminaries}
	\subsection{Basic definitions}
	For nonempty open sets $U,V\subset X$, let
	\[
	N_T(U,V)=\{n\in\N_0:T^nU\cap V\ne\varnothing\}.
	\]
	We omit the subscript when the map is clear. The system is \emph{transitive} if every such set is nonempty, \emph{minimal} if every forward orbit is dense, and \emph{weakly mixing} if $T\times T$ is transitive. Minimality implies surjectivity: the compact set $T(X)$ contains the forward orbit of $Tx$ for any $x\in X$, and hence equals $X$. Equivalently, mixing along an infinite set $A$ means that $A\setminus N_T(U,V)$ is finite for every nonempty open $U,V$.

	For a sequence $(a_j)_{j\ge1}$ of positive integers, write
	\[
	\operatorname{FS}(a_j)=\left\{\sum_{j\in F}a_j:
	\varnothing\ne F\subset\N\text{ finite}\right\}.
	\]
	An \emph{IP-set} is a set containing some $\operatorname{FS}(a_j)$; an \emph{IP$^*$-set} is a set that meets every IP-set.
	A minimal system is \emph{mildly mixing} if
	$N_T(U,V)$ is an IP$^*$-set for all nonempty
	open sets $U,V\subset X$. Equivalently, its
	product with every transitive system is
	transitive; see \cite{HY04,HSY04}.
	In the minimal setting, mild mixing is also equivalent to requiring every finite intersection of open-set return-time sets to be IP$^*$ \cite{HY04}.
	
	A set $S\subset\N_0$ is \emph{syndetic} if there is $L\ge0$ such that every interval $\{m,m+1,\ldots,m+L\}$ meets $S$, and \emph{thick} if it contains intervals of arbitrarily large finite length. We use the standard facts that $N_T(x,U)=\{n\ge0:T^nx\in U\}$ is syndetic in a minimal system, and that $N_T(U,V)$ is thick in a weakly mixing system; see \cite{Akin97}. In particular, a syndetic set meets every thick set.
	
	For a finite open cover $\alpha$, let $N(\alpha)$ be the least cardinality of a subcover, and define
	\[
	h(T,\alpha)=\lim_{n\to\infty}\frac1n
	\log N\left(\bigvee_{j=0}^{n-1}T^{-j}\alpha\right).
	\]
	The system has \emph{uniformly positive entropy} (UPE) if $h(T,\alpha)>0$ for every standard binary open cover $\alpha$, as defined in the introduction.
	
	A pair of distinct points $(x,y)$ is an \emph{entropy pair} if
	\[
	h(T,\{X\setminus K_x,X\setminus K_y\})>0
	\]
	for every pair of disjoint closed neighborhoods $K_x,K_y$ of $x,y$, respectively. Write $E(X,T)$ for the set of entropy pairs and $\Delta_X$ for the diagonal. We shall use the standard characterization \cite{B93,GW94}
	\begin{equation}\label{eq:upepairs}
		(X,T)\text{ has UPE}
		\quad\Longleftrightarrow\quad
		E(X,T)=X^2\setminus\Delta_X.
	\end{equation}
	For the forward asymptotic relation $\Asy(T)$ defined in the introduction,
	Blanchard, Host and Ruette proved, for homeomorphisms,
	\begin{equation}\label{eq:entropyasymptotic}
		E(X,T)\subseteq\overline{\Asy(T)};
	\end{equation}
	see \cite[Proposition~4]{BHR02}. This result will be used to prove Theorem~\ref{thm:upe}.
	
	We also use the following multiple recurrence property of minimal weakly mixing homeomorphisms \cite[Theorems~1 and~5]{KO12}: for every $r\ge1$ and nonempty open sets $U_0,\ldots,U_r\subset X$, there are $x\in X$ and $d\ge1$ such that
	\begin{equation}\label{eq:multiplehitting}
		x\in U_0,\qquad T^{id}x\in U_i\quad(1\le i\le r).
	\end{equation}
	Equivalently, $\bigcup_{d\ge1}\{(x,T^dx,\ldots,T^{rd}x):x\in X\}$ is dense in $X^{r+1}$.
	
	 \subsection{Natural extensions}
	
	We recall the following properties of natural extensions.
	The preservation of UPE was established by He and Liu
	\cite{HL99}. We include a proof of the distance property
	needed below.
	
	\begin{lem}\label{lem:naturalextension}
		Let $T:X\to X$ be a continuous surjection, and let
		\[
		\widehat X=\{(x_j)_{j\ge0}\in X^{\N_0}:Tx_{j+1}=x_j\},\qquad
		\widehat T(x_0,x_1,\ldots)=(Tx_0,x_0,x_1,\ldots).
		\]
		Let $\pi_m:\widehat X\to X$ be the $m$th coordinate map. Then $\widehat T$ is a homeomorphism, each $\pi_m$ is onto, and $\pi_m\widehat T=T\pi_m$. Moreover:
		\begin{enumerate}
			\item $(X,T)$ is minimal if and only if $(\widehat X,\widehat T)$ is minimal.
			\item $(X,T)$ has UPE if and only if $(\widehat X,\widehat T)$ has UPE.
			\item For every $m\ge0$ and nonempty open $U,V\subset X$,
			\begin{equation}\label{eq:extensionhitting}
				N_{\widehat T}(\pi_m^{-1}U,\pi_m^{-1}V)=N_T(U,V).
			\end{equation}
			Consequently, weak mixing and mixing along any fixed infinite set $A$ hold for $(X,T)$ if and only if they hold for its natural extension.
			\item Given nonempty open $\widehat U,\widehat V\subset\widehat X$ and $\varepsilon>0$, there are nonempty open $U,V\subset X$ and $\delta>0$ such that, for every $n\ge0$,
			\begin{equation}\label{eq:extensiondistance}
				\dist(T^nU,T^nV)<\delta
				\quad\Longrightarrow\quad
				\dist(\widehat T^n\widehat U,\widehat T^n\widehat V)<\varepsilon.
			\end{equation}
			In particular, the distance condition \eqref{eq:distance} along a fixed infinite set $A$ holds for $(X,T)$ if and only if it holds for $(\widehat X,\widehat T)$.
		\end{enumerate}
	\end{lem}
	
	\begin{proof}
		Assertions (1) and (3) are standard consequences
		of the cylinder basis of $\widehat X$, while (2)
		follows from \cite{HL99}. We prove (4).
		
		Equip $\widehat X$ with the compatible metric
		\[
		\widehat\rho(\hat x,\hat y)
		=\sum_{k=0}^{\infty}2^{-k-1}
		\min\{1,\rho(x_k,y_k)\}.
		\]
		Choose $r\ge0$ and nonempty open sets
		$U,V\subset X$ such that
		\[
		\pi_r^{-1}U\subset\widehat U,
		\qquad
		\pi_r^{-1}V\subset\widehat V,
		\]
		and $\sum_{k>r}2^{-k-1}<\frac{\varepsilon}{2}$.
		By uniform continuity, there is $\delta>0$
		such that
		\[
		\rho(a,b)<\delta
		\quad\Longrightarrow\quad
		\rho(T^ia,T^ib)<\frac{\varepsilon}{2}
		\qquad(0\le i\le r).
		\]
		
		Suppose that $\dist(T^nU,T^nV)<\delta$.
		Choose $u\in U$ and $v\in V$ with
		$\rho(T^nu,T^nv)<\delta$, and take
		$\hat u\in\pi_r^{-1}(u)$ and
		$\hat v\in\pi_r^{-1}(v)$.
		For $0\le k\le r$,
		\[
		\pi_k(\widehat T^n\hat u)=T^{n+r-k}u,
		\qquad
		\pi_k(\widehat T^n\hat v)=T^{n+r-k}v.
		\]
		Hence
		\[
		\widehat\rho(\widehat T^n\hat u,
		\widehat T^n\hat v)
		<
		\frac{\varepsilon}{2}
		+\sum_{k>r}2^{-k-1}
		<\varepsilon,
		\]
		which proves \eqref{eq:extensiondistance}.
		
		The converse implication for
		condition~\eqref{eq:distance}
		follows from the uniform continuity of $\pi_0$.
	\end{proof}
	
	\section{Mixing}\label{sec:criterion}
	
	\subsection{A criterion for mixing}
	We first work with homeomorphisms. Fix a countable base $\mathcal B$ of nonempty open subsets of $X$. For a strictly increasing sequence $(n_j)$ of nonnegative integers, consider the following hypothesis:
	\begin{equation*}\tag{\textbf{H1}}\label{hyp:H1}
		\{j\in\N:T^{n_j}U\cap V\ne\varnothing\}
		\text{ is finite or cofinite in }\N
		\qquad(U,V\in\mathcal B).
	\end{equation*}
	Every strictly increasing sequence has a subsequence satisfying \eqref{hyp:H1}. Indeed, enumerate $\mathcal B\times\mathcal B$ and successively choose an infinite subsequence along which the corresponding intersection is always nonempty or always empty. A diagonal subsequence has the required property. The threshold may depend on $U,V$.
	
	For a sequence satisfying \eqref{hyp:H1} and a point $x\in X$, define
	\begin{equation}\label{eq:limits}
		L(x)=\left\{y\in X:
		\begin{array}{l}
			\text{there exist }j_k\to\infty\text{ and }x_k\to x\text{ such that}\\
			T^{n_{j_k}}x_k\to y
		\end{array}\right\}.
	\end{equation}
	The points $x_k$ are allowed to vary. The following elementary properties will be used.
	
	\begin{lem}\label{lem:limits}
		Let $T$ be a homeomorphism and let $(n_j)$ satisfy \eqref{hyp:H1}. Each $L(x)$ is a nonempty compact set. Moreover:
		\begin{enumerate}
			\item If $x_k\to x$, $y_k\to y$ and $y_k\in L(x_k)$, then $y\in L(x)$.
			\item $L(T^rx)=T^rL(x)$ for every $r\in\Z$.
			\item If $O,V\subset X$ are open and $L(x)\cap V\ne\varnothing$ for some $x\in O$, then $T^{n_j}O\cap V\ne\varnothing$ for all sufficiently large $j$.
			\item There is a dense $G_\delta$ set $D\subset X$ such that, whenever $x\in D$ and an open set $V$ meets $L(x)$, some neighborhood $O$ of $x$ satisfies
			\begin{equation}\label{eq:lower}
				L(v)\cap V\ne\varnothing\qquad(v\in O).
			\end{equation}
		\end{enumerate}
	\end{lem}
	
	\begin{proof}
		Nonemptiness follows by taking an accumulation point of $(T^{n_j}x)$. To prove the first assertion, for each $k$ choose $j_k\ge k$ and $u_k\in X$ such that
		\[
		\rho(u_k,x_k)<1/k,\qquad
		\rho(T^{n_{j_k}}u_k,y_k)<1/k.
		\]
		Then $u_k\to x$ and $T^{n_{j_k}}u_k\to y$, so $y\in L(x)$. This also proves that $L(x)$ is closed, hence compact. The second assertion follows from \eqref{eq:limits} and the commutation of the iterates of $T$; applying the same argument to $T^{-r}$ gives equality.
		
		For the third assertion, choose $y\in L(x)\cap V$ and basic neighborhoods $U_0,V_0$ with $x\in U_0\subset O$ and $y\in V_0\subset V$. Definition~\eqref{eq:limits} gives $T^{n_j}U_0\cap V_0\ne\varnothing$ for infinitely many $j$. By \eqref{hyp:H1}, this holds for all sufficiently large $j$.

		For the last assertion, choose a countable family
		$\mathcal C$ of closed balls such that, for every
		$x\in X$ and every open neighborhood $U$ of $x$,
		there exists $K\in\mathcal C$ satisfying
		$x\in\operatorname{int}K\subset K\subset U.$
		For $K\in\mathcal C$, the set
		\[
		A_K=\{x\in X:L(x)\cap K\ne\varnothing\}
		\]
		is closed by the first assertion and compactness. Hence
		\[
		D=X\setminus\bigcup_{K\in\mathcal C}\partial A_K
		\]
		is a dense $G_\delta$ set. If $x\in D$ and $y\in L(x)\cap V$, choose $K\in\mathcal C$ with $y\in\operatorname{int}K\subset K\subset V$. Then $x\in A_K\setminus\partial A_K=\operatorname{int}A_K$, which gives \eqref{eq:lower}.
	\end{proof}
	
	The next lemma will be used to pass from $s$-fold to $(s+1)$-fold intersections.
	
	\begin{lem}\label{lem:completion}
		Let $T$ be a homeomorphism, let $(n_j)$ satisfy \eqref{hyp:H1}, and let $L$ and $D$ be as in Lemma~\ref{lem:limits}. Suppose that $s\ge2$ and
		\begin{equation}\label{eq:sintersection}
			\bigcap_{i=1}^sL(v_i)\ne\varnothing
			\qquad\text{for every }(v_1,\ldots,v_s)\in X^s.
		\end{equation}
		If $x_1,\ldots,x_s\in D$, there is $z\in\bigcap_{i=1}^sL(x_i)$ such that
		\begin{equation}\label{eq:outputfirst}
			L(z)\cap\bigcap_{i=1}^sL(y_i)\ne\varnothing
			\qquad\text{for every }(y_1,\ldots,y_s)\in\prod_{i=1}^sL(x_i).
		\end{equation}
	\end{lem}
	
	\begin{proof}
		Put $K=\bigcap_{i=1}^sL(x_i)$, which is nonempty and compact by \eqref{eq:sintersection} and Lemma~\ref{lem:limits}. For $\mathbf y=(y_1,\ldots,y_s)\in\prod_iL(x_i)$, let
		\[
		K_{\mathbf y}=\left\{z\in K:
		L(z)\cap\bigcap_{i=1}^sL(y_i)\ne\varnothing\right\}.
		\]
		The set $\{(z,w)\in K\times X:w\in L(z)\cap\bigcap_iL(y_i)\}$ is compact by Lemma~\ref{lem:limits}(1). Its projection onto the first coordinate is $K_{\mathbf y}$, so $K_{\mathbf y}$ is closed. We prove that these sets have the finite intersection property.
		
		Fix $q\ge1$ and points $y_i^\ell\in L(x_i)$, where $1\le i\le s$ and $1\le\ell\le q$. For each positive integer $m$, choose open neighborhoods
		\[
		y_i^\ell\in V_{i,m}^\ell\subset B_{1/m}(y_i^\ell).
		\]
		By \eqref{eq:lower}, there are open neighborhoods
		\[
		x_i\in O_{i,m}^0\subset B_{1/m}(x_i)
		\]
		such that $L(v)\cap V_{i,m}^\ell\ne\varnothing$ for every $v\in O_{i,m}^0$ and every $\ell$.
		
		For $1\le\ell\le q$, suppose that nonempty open sets $O_{i,m}^{\ell-1}\subset O_{i,m}^0$ have been chosen. Every $v\in O_{i,m}^{\ell-1}$ satisfies $L(v)\cap V_{i,m}^\ell\ne\varnothing$. Thus Lemma~\ref{lem:limits}(3) gives, for all sufficiently large $j$,
		\[
		T^{n_j}O_{i,m}^{\ell-1}\cap V_{i,m}^\ell\ne\varnothing
		\qquad(1\le i\le s).
		\]
		Choose sufficiently large $j_m^\ell\ge m$ such that this holds for $1\le i\le s$, and put
		\begin{equation}\label{eq:nested}
			O_{i,m}^\ell=O_{i,m}^{\ell-1}\cap T^{-n_{j_m^\ell}}V_{i,m}^\ell.
		\end{equation}
		This defines nonempty open sets $O_{i,m}^\ell$ for $1\le\ell\le q$, with $O_{i,m}^q\subset O_{i,m}^\ell\subset O_{i,m}^0$.
		
		Choose $v_{i,m}\in O_{i,m}^q$ and, by \eqref{eq:sintersection}, choose
		\[
		z_m\in\bigcap_{i=1}^sL(v_{i,m}).
		\]
		For $1\le\ell\le q$, set
		\[
		w_m^\ell=T^{n_{j_m^\ell}}z_m,
		\qquad
		u_{i,m}^\ell=T^{n_{j_m^\ell}}v_{i,m}.
		\]
		Equivariance and \eqref{eq:nested} give
		\begin{equation}\label{eq:finitewitnesses}
			w_m^\ell\in L(u_{i,m}^\ell),
			\qquad u_{i,m}^\ell\in V_{i,m}^\ell.
		\end{equation}
		By compactness of $X^{q+1}$, there is a subsequence $m_k\to\infty$ such that
		\[
		z_{m_k}\to z,\qquad
		T^{n_{j_{m_k}^\ell}}z_{m_k}=w_{m_k}^\ell\to w^\ell
		\quad(1\le\ell\le q).
		\]
		Since $j_{m_k}^\ell\ge m_k\to\infty$, definition~\eqref{eq:limits} gives $w^\ell\in L(z)$ for every $\ell$. The neighborhood choices give $v_{i,m_k}\to x_i$ and $u_{i,m_k}^\ell\to y_i^\ell$. Lemma~\ref{lem:limits}(1), applied to $z_{m_k}\in L(v_{i,m_k})$ and \eqref{eq:finitewitnesses}, yields $z\in K$ and $w^\ell\in\bigcap_iL(y_i^\ell)$. Hence
		\[
		z\in\bigcap_{\ell=1}^qK_{(y_1^\ell,\ldots,y_s^\ell)}.
		\]
		Since the finite family of tuples was arbitrary, compactness gives $\bigcap_{\mathbf y\in\prod_iL(x_i)}K_{\mathbf y}\ne\varnothing$. Any $z$ in this intersection satisfies \eqref{eq:outputfirst} for every $\mathbf y\in\prod_iL(x_i)$.
	\end{proof}
	\begin{prop}\label{prop:alllimits}
		Let $(X,T)$ be a minimal system with $T$ a homeomorphism, and let $(n_j)$ satisfy \eqref{hyp:H1}. Suppose that
		\begin{equation}\label{eq:pairintersection}
			L(x)\cap L(y)\ne\varnothing\qquad(x,y\in X).
		\end{equation}
		Then $L(x)=X$ for every $x\in X$.
	\end{prop}
	
	\begin{proof}
		First we show that $(X,T)$ is weakly mixing. Given nonempty open sets $U_1,U_2,V_1,V_2$, choose $k\ge0$ such that $W=V_1\cap T^{-k}V_2\ne\varnothing$. Choose $a\in U_1$, $b\in T^{-k}U_2$ and $z\in L(a)\cap L(b)$. By minimality and continuity, there are $t\ge0$ and an open neighborhood $O$ of $z$ such that $T^tO\subset W$. Lemma~\ref{lem:limits}(3) gives, for all sufficiently large $j$,
		\[
		T^{n_j}U_1\cap O\ne\varnothing,\qquad
		T^{n_j}(T^{-k}U_2)\cap O\ne\varnothing.
		\]
		Applying $T^t$ to the first intersection and $T^{t+k}$ to the second gives simultaneous hits from $U_1,U_2$ to $V_1,V_2$ at time $n_j+t$. Thus $T\times T$ is transitive.
		
		We prove \eqref{eq:sintersection} for all $s\ge2$ by induction. The case $s=2$ is the hypothesis. Suppose that it holds for some $s\ge2$, and put
		\[
		E_{s+1}=\left\{(v_0,\ldots,v_s)\in X^{s+1}:\bigcap_{i=0}^sL(v_i)\ne\varnothing\right\}.
		\]
		This set is closed by Lemma~\ref{lem:limits}(1) and compactness.
		
		Fix $d\ge1$. Since $T$ is a homeomorphism and $D$ is a dense $G_\delta$ set, we can choose
		\[
		b\in\bigcap_{i=0}^{s-1}T^{-id}D.
		\]
		Apply Lemma~\ref{lem:completion} to $b,T^db,\ldots,T^{(s-1)d}b$. It gives a point
		\[
		z\in\bigcap_{i=0}^{s-1}L(T^{id}b)
		\]
		for which \eqref{eq:outputfirst} holds. Now let
		\[
		y_0=T^{-d}z,\qquad y_i=T^{id}z\quad(1\le i<s).
		\]
		Since $z\in L(T^db)=T^dL(b),$
		we have $y_0\in L(b)$ and $y_i\in L(T^{id}b)$ for $1\le i<s$. By \eqref{eq:outputfirst},
		\[
		L(z)\cap L(T^{-d}z)\cap L(T^dz)\cap\cdots\cap L(T^{(s-1)d}z)
		\ne\varnothing.
		\]
		Consequently
		\[
		A_d=\left\{x\in X:\bigcap_{i=0}^sL(T^{id}x)\ne\varnothing\right\}
		\]
		contains $T^{-d}z$. It is closed because it is the inverse image of $E_{s+1}$ under the continuous map $x\mapsto(x,T^dx,\ldots,T^{sd}x)$. Equivariance gives $TA_d=A_d$, and minimality therefore gives $A_d=X$.
		
		Since $d\ge1$ was arbitrary, $E_{s+1}$ contains
		\[
		\bigcup_{d\ge1}\{(x,T^dx,\ldots,T^{sd}x):x\in X\}.
		\]
		This union is dense in $X^{s+1}$ by \eqref{eq:multiplehitting}. Closedness of $E_{s+1}$ now gives $E_{s+1}=X^{s+1}$, proving the induction step.
		
		The family $\{L(x):x\in X\}$ therefore has the finite intersection property, so
		\[
		Z=\bigcap_{x\in X}L(x)\ne\varnothing.
		\]
		This is a closed set, and equivariance gives $TZ=Z$. Minimality yields $Z=X$, and therefore $L(x)=X$ for every $x$.
	\end{proof}

	\subsection{Applications}\label{sec:applications}

	\begin{proof}[\textbf{Proof of  $(1)\Leftrightarrow(2)$ of Theorem~\ref{thm:criterion}}]
		$(1)\Rightarrow(2)$.
		For any nonempty open sets $U,V$ and $\varepsilon>0$, choose a
		nonempty open set $W$ with $\diam W<\varepsilon$.
		For all sufficiently large $n\in A$, both $T^nU$ and $T^nV$ meet $W$,
		which proves \eqref{eq:distance}.
		
		$(2)\Rightarrow(1)$.
		First suppose that $T$ is a homeomorphism. If mixing along $A$ fails,
		there are nonempty open sets $U,V$ and a strictly increasing sequence
		$(n_j)$ in $A$ such that
		\begin{equation}\label{eq:badtimes}
			T^{n_j}U\cap V=\varnothing\qquad(j\ge1).
		\end{equation}
		Pass to a subsequence satisfying \eqref{hyp:H1}, and define $L$ by
		\eqref{eq:limits}.
		
		For any $x,y\in X$, condition~(2) gives $j_m\ge m$ and points
		$a_m\in B_{1/m}(x)$, $b_m\in B_{1/m}(y)$ with
		\[
		\rho(T^{n_{j_m}}a_m,T^{n_{j_m}}b_m)<1/m.
		\]
		By compactness, the two image sequences have a common subsequential
		limit in $L(x)\cap L(y)$. Proposition~\ref{prop:alllimits} gives
		$L(x)=X$ for every $x$. Lemma~\ref{lem:limits}(3), applied to a point
		$x\in U$ and the open set $V$, now contradicts \eqref{eq:badtimes}.
		
		For a general minimal system, Lemma~\ref{lem:naturalextension} shows
		that the natural extension is minimal and satisfies the same distance
		condition along $A$. The homeomorphism case gives mixing along $A$
		for the extension, and hence for $(X,T)$.
	\end{proof}

	\begin{cor}\label{thm:asy}
		Let $(X,T)$ be minimal. If $\Asy(T)$ is dense in $X^2$,
		then $(X,T)$ is mixing. For a homeomorphism, the same
		conclusion holds if $\Asy(T^{-1})$ is dense.
	\end{cor}
	
	\begin{proof}
		If $\Asy(T)$ is dense in $X^2$, then for any nonempty
		open sets $U,V$, there exist $x\in U$ and $y\in V$
		such that
		\[
		\rho(T^nx,T^ny)\longrightarrow 0.
		\]
		Hence
		\[
		\operatorname{dist}(T^nU,T^nV)\longrightarrow 0.
		\]
		The conclusion follows from Theorem~\ref{thm:criterion}.
		The case of $\Asy(T^{-1})$ follows by applying the
		same argument to $T^{-1}$.
	\end{proof}
	\begin{proof}[\textbf{Proof of Theorem~\ref{thm:upe}}]
		By Lemma~\ref{lem:naturalextension}, it suffices to treat a minimal homeomorphism with UPE. Equations~\eqref{eq:upepairs} and~\eqref{eq:entropyasymptotic} give $\overline{\Asy(T)}=X^2$. Corollary~\ref{thm:asy} gives mixing, which then passes to the original system.
	\end{proof}
	
	Following \cite{KSS03,GY09}, a system is \emph{KSS} if there is a closed equivalence relation $R\subset X^2$ such that
	\begin{equation}\label{eq:kss}
		(T\times T)R\subset R,\qquad
		\bigcap_{n\ge0}(T\times T)^nR=\Delta_X,\qquad
		\overline{\bigcup_{n\ge0}(T\times T)^{-n}R}=X^2.
	\end{equation}
	This is the Kolmogorov property defined through closed equivalence relations in \cite{KSS03}; it is different from the all-orders UPE property mentioned in the introduction.
	
	Glasner and Ye \cite[Theorem~13.8(1)]{GY09}
	proved that every minimal KSS system is weakly mixing.
	The following result strengthens this conclusion to mixing.
	
	\begin{cor}\label{cor:kss}
		Every minimal KSS system is topologically mixing.
	\end{cor}
	
	\begin{proof}
		The first two conditions in \eqref{eq:kss} imply $R\subset\Asy(T)$.  If $(T\times T)^n(x,y)\in R$ for some $n\ge0$, then $(x,y)\in\Asy(T)$. Thus the third condition gives $\overline{\Asy(T)}=X^2$. Apply Corollary~\ref{thm:asy}.
	\end{proof}

	A minimal homeomorphism is \emph{doubly minimal} if, whenever $x,y$ belong to different $T$-orbits, their two-sided orbit under $T\times T$ is dense in $X^2$.
	
	\begin{cor}\label{cor:double}
		Every infinite doubly minimal homeomorphism is topologically mixing.
	\end{cor}
	
	\begin{proof}
		Fix nonempty open sets $U,V$ and $\varepsilon>0$. King's boundedness theorem \cite[Theorem~2.6]{King90} gives points $x,y$ in different $T$-orbits such that
		\[
		\rho(T^mx,T^my)\le\varepsilon\qquad(m\ge1).
		\]
		By double minimality there is $k\in\Z$ with $T^kx\in U$ and $T^ky\in V$. Thus
		\[
		\dist(T^nU,T^nV)
		\le \rho(T^{n+k}x,T^{n+k}y)
		\le\varepsilon\qquad(n+k\ge1).
		\]
		Since $\varepsilon$ was arbitrary, condition~\eqref{eq:distance} holds with $A=\N_0$. Theorem~\ref{thm:criterion} gives mixing.
	\end{proof}

	\section{Binary-cover complexity}\label{sec:binary}
	
	We use the cover complexity $c_\alpha(F)$ and the difference set $\Delta(B)$ defined in the introduction.
	
	\subsection{Binary obstructions}
	
	The following lemma gives a common obstruction to unbounded binary-cover complexity and to $\Delta^*$-transitivity.
	
	\begin{lem}\label{lem:binarynames}
		Let $(X,T)$ be minimal and let $A\subset\N_0$ be infinite. If $(X,T)$ is not mixing along $A$, then there are disjoint closed neighborhoods $K_0,K_1\subset X$, an infinite $B\subset A$, and an open cover $\{Q_0,Q_1\}$ of $X$ such that, with $A_i=X\setminus K_i$,
		\begin{equation}\label{eq:fixednames}
			T^bQ_i\subset A_i\qquad(b\in B,\ i=0,1).
		\end{equation}
		The standard binary open cover $\alpha=\{A_0,A_1\}$ satisfies
		\begin{equation}\label{eq:binarytwo}
			c_\alpha(F)=2\qquad\text{for every nonempty finite }F\subset B.
		\end{equation}
		Moreover, writing $U_i=\operatorname{int}K_i$, we have
		\begin{equation}\label{eq:differencegap}
			\Delta(B)\cap\bigl(N(U_0,U_1)\cup N(U_1,U_0)\bigr)=\varnothing.
		\end{equation}
	\end{lem}
	
	\begin{proof}
		First suppose that $T$ is a homeomorphism. Mixing along $A$ is equivalent for $T$ and $T^{-1}$.   The equivalence of (1) and (2) in Theorem~\ref{thm:criterion},
		already proved in Section~\ref{sec:criterion}, applied to $T^{-1}$, gives nonempty open $U,V$, an $\varepsilon>0$ and an infinite $S\subset A$ such that
		\begin{equation}\label{eq:backseparation}
			\dist(T^{-n}U,T^{-n}V)\ge\varepsilon\qquad(n\in S).
		\end{equation}
		In particular, $U\cap V=\varnothing$. Choose closed neighborhoods $K_0\subset U$, $K_1\subset V$ with nonempty interiors, and a finite open cover $\{W_1,\ldots,W_M\}$ with $\diam W_j<\varepsilon$. For each $n\in S$, the set $W_j$ cannot meet both $T^{-n}K_0$ and $T^{-n}K_1$. We can therefore choose $e_j(n)\in\{0,1\}$ such that
		\[
		W_j\subset T^{-n}A_{e_j(n)}.
		\]
		There are only finitely many possible tuples $(e_1(n),\ldots,e_M(n))$.
		Choose an infinite $B\subset S$ on which the tuple is constant,
		say equal to $(e_1,\ldots,e_M)$, and put
		\[
		Q_i=\bigcup_{e_j=i}W_j\qquad(i=0,1).
		\]
		Then $Q_0\cup Q_1=X$, and \eqref{eq:fixednames} follows.
		
		Now suppose only that $T$ is continuous. Its natural extension is minimal and is not mixing along $A$, by Lemma~\ref{lem:naturalextension}. Apply the preceding construction to obtain a standard binary open cover $\{\widehat A_0,\widehat A_1\}$, an open cover $\{\widehat Q_0,\widehat Q_1\}$ and an infinite $B\subset A$ such that
		\[
		\widehat T^b\widehat Q_i\subset\widehat A_i
		\qquad(b\in B,\ i=0,1).
		\]
		Since neither $\widehat A_i$ is dense, the cylinder basis gives a common coordinate $m$ and nonempty open sets $V_0,V_1\subset X$ such that
		\[
		\pi_m^{-1}V_i\subset\widehat X\setminus\overline{\widehat A_i}
		\qquad(i=0,1).
		\]
		The equality $\widehat A_0\cup\widehat A_1=\widehat X$ and
		surjectivity of $\pi_m$ give $V_0\cap V_1=\varnothing$. Choose closed neighborhoods $L_i\subset V_i$ with nonempty interiors. By compactness, there is a finite cylinder cover $\{\pi_r^{-1}W_j:1\le j\le M\}$, at a common coordinate $r\ge m$, refining $\{\widehat Q_0,\widehat Q_1\}$. Choose $e_j\in\{0,1\}$ with $\pi_r^{-1}W_j\subset\widehat Q_{e_j}$ and put $Q_i=\bigcup_{e_j=i}W_j$. Then
		\[
		\pi_r^{-1}Q_i\subset\widehat Q_i,\qquad Q_0\cup Q_1=X.
		\]
		Put $K_i=T^{-(r-m)}L_i$ and $A_i=X\setminus K_i$. The sets $K_i$ are disjoint closed neighborhoods with nonempty interiors, because $T^{r-m}$ is onto. The identity $\pi_m=T^{r-m}\pi_r$ gives
		\[
		\widehat A_i\subset\pi_m^{-1}(X\setminus L_i)=\pi_r^{-1}A_i.
		\]
		Consequently, for every $b\in B$,
		\[
		\pi_r^{-1}Q_i\subset\widehat Q_i
		\subset\widehat T^{-b}\widehat A_i
		\subset\widehat T^{-b}\pi_r^{-1}A_i
		=\pi_r^{-1}(T^{-b}A_i).
		\]
		Surjectivity of $\pi_r$ gives $Q_i\subset T^{-b}A_i$ for every $b\in B$, proving \eqref{eq:fixednames} with the same set $B$.
		
		In either case, \eqref{eq:fixednames} gives
		\begin{equation}\label{eq:twoconstantnames}
			X=\left(\bigcap_{b\in B}T^{-b}A_0\right)
			\cup\left(\bigcap_{b\in B}T^{-b}A_1\right).
		\end{equation}
		Thus $c_\alpha(F)\le2$ for every nonempty finite
		$F\subset B$. Since $T$ is onto and neither $A_i$
		equals $X$, no member of
		$\bigvee_{b\in F}T^{-b}\alpha$ covers $X$.
		Hence $c_\alpha(F)=2$.
		
		Finally, if $b<c$ belong to $B$ and
		$T^{c-b}U_0\cap U_1\ne\varnothing$, then
		$T^{-b}K_0\cap T^{-c}K_1\ne\varnothing$
		by surjectivity, contradicting
		\eqref{eq:twoconstantnames}. Interchanging
		$0$ and $1$ proves \eqref{eq:differencegap}.
	\end{proof}

	\subsection{Proof of Theorem~\ref{thm:criterion}}\label{sec:binaryproof}

	\begin{proof}
		$(1)\Rightarrow(3)$.
		Fix a standard binary open cover $\alpha=\{A_0,A_1\}$ and an infinite
		$S\subset A$. Choose nonempty open sets
		\[
		V_i\subset X\setminus\overline{A_i}\qquad(i=0,1).
		\]
		By surjectivity and mixing along $A$, we can inductively choose
		$t_1<t_2<\cdots$ in $S$ such that
		\begin{equation}\label{eq:binaryindependence}
			C_\omega=\bigcap_{j=1}^kT^{-t_j}V_{\omega_j}\ne\varnothing
			\qquad(\omega\in\{0,1\}^k,\ k\ge1).
		\end{equation}
		Indeed, after $t_1,\ldots,t_k$ have been chosen,
		$T^nC_\omega$ meets both $V_0$ and $V_1$ for all sufficiently large
		$n\in A$, simultaneously for every $\omega$.
		
		Since $V_i\cap A_i=\varnothing$, each member of
		$\bigvee_{j=1}^kT^{-t_j}\alpha$ meets at most one $C_\omega$. Hence
		\[
		c_\alpha(\{t_1,\ldots,t_k\})=2^k.
		\]
		Writing $S=\{s_1<s_2<\cdots\}$, monotonicity under refinement gives
		\[
		c_\alpha(\{s_1,\ldots,s_m\})\ge2^k
		\qquad\text{whenever }s_m\ge t_k,
		\]
		which proves (3). Restricting \eqref{eq:binaryindependence} to any
		nonempty finite $F\subset\{t_j:j\ge1\}$ gives
		$c_\alpha(F)=2^{|F|}$, as stated in
		Remark~\ref{rem:mixingconsequences}.
		
		$(3)\Rightarrow(1)$.
		This follows from Lemma~\ref{lem:binarynames} and
		\eqref{eq:binarytwo}.
		
		Now let $A=\N_0$.
		The implication $(4)\Rightarrow(1)$ follows from
		Lemma~\ref{lem:binarynames} and \eqref{eq:differencegap}.
		For $(1)\Rightarrow(5)$, mixing passes to the Cartesian square,
		and every cofinite set meets $\Delta(B)$ for every infinite $B$.
		Finally, $(5)\Rightarrow(4)$ follows by testing the open sets
		$U\times X$ and $V\times X$.
		
		For arbitrary infinite $A$, the assertion \eqref{eq:introbinary}
		follows from Lemma~\ref{lem:binarynames}.
	\end{proof}
	
	\section{Local mixing}\label{sec:local}
	
	For a minimal system $(X,T)$, consider the regionally proximal relation
	associated with the cofinite family \cite{HSY04},
	\begin{equation}\label{eq:localrelation}
		\mathcal C(T)=\left\{(x,y)\in X^2:
		\lim_{n\to\infty}\dist(T^nU,T^nV)=0
		\text{ for all open }U\ni x,\ V\ni y\right\}.
	\end{equation}
	Write $\mathcal C(T)[x]=\{y:(x,y)\in\mathcal C(T)\}$.
	Let $\pi_{\mathrm{eq}}:X\to X_{\mathrm{eq}}$ be the maximal equicontinuous factor map, and write
	\[
	R_{\pi_{\mathrm{eq}}}
	=\{(x,y)\in X^2:\pi_{\mathrm{eq}}(x)=\pi_{\mathrm{eq}}(y)\}.
	\]

	\subsection{Local relations and cyclic factors}
	
	\begin{lem}\label{lem:localproperties}
		The relation $\mathcal C(T)$ is closed and forward invariant under $T\times T$.
		For every integer $q\ge1$,
		\[
		\mathcal C(T^q)=\mathcal C(T)\subset R_{\pi_{\mathrm{eq}}}.
		\]
	\end{lem}
	
	\begin{proof}
		Closedness is immediate. Forward invariance follows by applying
		the defining condition to $T^{-1}U$ and $T^{-1}V$ for neighborhoods
		$U,V$ of $Tx,Ty$, respectively, and using uniform continuity.
		The equality $\mathcal C(T^q)=\mathcal C(T)$ follows by restricting
		to multiples of $q$ and using uniform continuity of $T^r$,
		$0\le r<q$.
		
		Finally, distinct points in the maximal equicontinuous factor
		remain uniformly separated under common iterates in a compatible
		invariant metric. Uniform continuity of $\pi_{\mathrm{eq}}$
		therefore gives $\mathcal C(T)\subset R_{\pi_{\mathrm{eq}}}$.
	\end{proof}
	
	We also use the following consequence of the maximal
	equicontinuous factor; compare
	\cite[Propositions~29 and~30]{BHS08}.

	\begin{lem}\label{lem:finitefactor}
		Let $(X,T)$ be minimal and suppose that its maximal equicontinuous
		factor is finite. Then this factor is a $p$-cycle, and its fibers
		give a clopen decomposition
		\[
		X=\bigsqcup_{i\in\Z/p\Z}C_i,\qquad TC_i=C_{i+1}.
		\]
		Each map $S_i=T^p|_{C_i}$ is minimal and weakly mixing.
	\end{lem}
	
	\begin{proof}
		The finite factor is a cycle by minimality. Each $S_i$ is minimal,
		since the return times of any $x\in C_i$ to $C_i$ are precisely
		the multiples of $p$.
		
		Suppose that $S_0$ admits a nontrivial equicontinuous factor
		$\phi:(C_0,S_0)\to(Y,R)$. Define
		\[
		\widetilde R(y,i)=
		\begin{cases}
			(y,i+1),&i<p-1,\\
			(Ry,0),&i=p-1,
		\end{cases}
		\]
		and
		\[
		\widetilde\phi(T^iz)=(\phi(z),i),
		\qquad z\in C_0,\quad 0\le i<p.
		\]
		This is well defined: $T^iz=T^iz'$ implies
		$R\phi(z)=R\phi(z')$, hence $\phi(z)=\phi(z')$
		since $R$ is invertible. Each $T^i:C_0\to C_i$ is a continuous
		surjection between compact Hausdorff spaces, hence a quotient map.
		It follows that $\widetilde\phi$ is continuous.
		Thus $\widetilde\phi$ defines a nontrivial equicontinuous
		extension of the finite factor, contradicting its maximality.
		
		Hence $S_0$ has trivial maximal equicontinuous factor and is
		weakly mixing by \cite[Proposition~29]{BHS08}.
		Each $S_i$ is a factor of $S_0$ under $T^i:C_0\to C_i$,
		and is therefore weakly mixing.
	\end{proof}

	 \subsection{Local mixing and its consequences}\label{sec:localproof}
	
	\begin{proof}[\textbf{Proof of equivalence of (1)--(3) of  Theorem~\ref{thm:local}}]
		$(1)\Rightarrow(2)$.
		Choose $x\in O$ and, for a neighborhood $U$ of $x$, apply (1) to
		$U\cap O$ and $V$.
		
		$(2)\Rightarrow(3)$.
		Condition~(2) is equivalent to $O\subset\mathcal C(T)[x]$.
		By Lemma~\ref{lem:localproperties}, $x$ and $O$ lie in the same
		fiber of $\pi_{\mathrm{eq}}$. Minimality gives
		\[
		X=\bigcup_{n=0}^N T^{-n}O
		\]
		for some $N$. The images of these inverse images under
		$\pi_{\mathrm{eq}}$ are singletons, since the transformation on the factor is
		invertible. Thus $X_{\mathrm{eq}}$ is finite.
		Let $C_i$ and $S_i$ be as in Lemma~\ref{lem:finitefactor}, with $x\in C_0$.
		
		For nonempty open $U,V\subset C_0$, the set $N_{S_0}(x,U)$ is
		syndetic by minimality, whereas $N_{S_0}(O,V)$ is thick by weak mixing.
		They intersect, so
		\[
		\overline{\bigcup_{n\ge0}
			(S_0\times S_0)^n(\{x\}\times O)}=C_0^2.
		\]
		Closedness and forward invariance give $C_0^2\subset\mathcal C(T)$.
		Applying $T^i\times T^i$ and Lemma~\ref{lem:localproperties}, we obtain
		\begin{equation}\label{eq:localequality}
			\mathcal C(T)=R_{\pi_{\mathrm{eq}}}
			=\bigsqcup_{i\in\Z/p\Z}C_i^2.
		\end{equation}
		Hence $\mathcal C(S_i)=C_i^2$, and the distance criterion in
		Theorem~\ref{thm:criterion} makes every $S_i$ mixing.
		
		$(3)\Rightarrow(1)$.
		Mixing of each $S_i$ and Lemma~\ref{lem:localproperties} give
		\[
		\bigsqcup_i C_i^2\subset\mathcal C(T^p)
		=\mathcal C(T)\subset R_{\pi_{\mathrm{eq}}}
		\subset\bigsqcup_i C_i^2,
		\]
		where the last inclusion follows because the cyclic decomposition
		defines an equicontinuous factor. Taking $O=C_i$ proves (1).
		The equality identifies the components with the fibers of
		$\pi_{\mathrm{eq}}$, proving uniqueness up to cyclic relabeling.
		A nontrivial finite cyclic factor is incompatible with weak mixing,
		total minimality, or connectedness, so each of these assumptions forces
		$p=1$.
	\end{proof}
	
	The same proof shows that these conditions are equivalent to
	$\mathcal C(T)$ having nonempty interior: a nonempty open rectangle
	in $\mathcal C(T)$ gives a section with nonempty interior, while
	\eqref{eq:localequality} proves the converse. The equivalence with
	conditions (4) and (5), and the mildly mixing variant, will be proved
	in Section~\ref{sec:families}.
	
	\begin{cor}\label{cor:residues}
		Under the equivalent conditions of Theorem~\ref{thm:local},
		let $U,V\subset X$ be nonempty open sets and put
		\[
		I(U)=\{i\in\Z/p\Z:U\cap C_i\ne\varnothing\},
		\qquad D(U,V)=I(V)-I(U).
		\]
		For every sufficiently large $n$,
		\begin{equation}\label{eq:residues}
			T^nU\cap V\ne\varnothing
			\quad\Longleftrightarrow\quad n\bmod p\in D(U,V).
		\end{equation}
	\end{cor}
	
	\begin{proof}
		The forward implication follows from $TC_i=C_{i+1}$.
		Conversely, choose $0\le r<p$ and $i\in I(U)$, $j\in I(V)$
		with $r\equiv j-i\pmod p$. Surjectivity
		of $T^r:C_i\to C_j$ and mixing of $T^p|_{C_i}$ give
		\[
		T^{pk+r}U\cap V\ne\varnothing
		\]
		for all sufficiently large $k$. Since there are only
		finitely many residue classes, the result follows.
	\end{proof}
	
	\begin{cor}\label{cor:localasy}
		Let $(X,T)$ be minimal. If $\overline{\Asy(T)}$ has
		nonempty interior in $X^2$, then the decomposition
		in Theorem~\ref{thm:local} exists and
		\[
		\overline{\Asy(T)}
		=\bigsqcup_{i\in\Z/p\Z}C_i^2.
		\]
		It is sufficient that the closure of one stable class
		\[
		W^s(x)=\{y\in X:\rho(T^nx,T^ny)\longrightarrow0\}
		\]
		have nonempty interior. If $T$ is totally minimal,
		or if $X$ is connected, then $p=1$ and $T$ is mixing.
	\end{cor}
	
	\begin{proof}
		Put $F=\overline{\Asy(T)}\subset\mathcal C(T)$.
		If $F$ has nonempty interior, the equivalence of (1)--(3) in Theorem~\ref{thm:local}
		gives the cyclic decomposition. Since $F$ is closed
		and forward invariant, and contains a nonempty open
		rectangle in some $C_i^2$, transitivity of
		$S_i\times S_i$ gives $C_i^2\subset F$.
		Applying $T^r\times T^r$ yields
		\[
		F=\bigsqcup_{i\in\Z/p\Z}C_i^2.
		\]
		
		Since $W^s(x)\times W^s(x)\subset\Asy(T)$,
		nonempty interior of $\overline{W^s(x)}$ implies
		nonempty interior of $F$. Finally, total minimality
		or connectedness forces $p=1$.
	\end{proof}
	
	For homeomorphisms, the same conclusions hold for
	$\Asy(T^{-1})$ by applying Corollary~\ref{cor:localasy}
	to $T^{-1}$.

	\section{Proximality}\label{sec:families}\label{sec:onecell}
	
	We use proximal cells along infinite sets of times as in
	\cite{HSY04}. For nonempty open sets $U,V\subset X$ and
	$\varepsilon>0$, write
	\[
	C_\varepsilon(U,V)
	=\{n\ge0:\dist(T^nU,T^nV)<\varepsilon\}.
	\]
	\subsection{Furstenberg families}

	A \emph{proper Furstenberg family} $\mathcal F$ is a nonempty
	upward hereditary family of subsets of $\N_0$ not containing
	the empty set. Its dual is
	\[
	k\mathcal F
	=\{S\subset\N_0:S\cap F\ne\varnothing
	\text{ for every }F\in\mathcal F\}.
	\]
	The system is \emph{$\mathcal F$-transitive} if
	$N(U,V)\in\mathcal F$ for all nonempty open $U,V$,
	and \emph{$\mathcal F$-mixing} if $(X^2,T\times T)$ is
	$\mathcal F$-transitive. Recall that $P_S(x)$ is defined in \eqref{eq:introcell}
and that $M_T(X)$ denotes the set of $T$-invariant
Borel probability measures on $X$.
	
	For an infinite $S\subset\N_0$, let $RP_S$ consist of the
	pairs $(x,y)$ such that, for every pair of neighborhoods
	$U,V$ of $x,y$ and every $\varepsilon>0$, $	C_\varepsilon(U,V)\cap S\ne\varnothing.$
	This relation is closed and forward invariant. We write
	$RP_S[x]=\{y:(x,y)\in RP_S\}$.
	When $\mathcal F$ contains the family
	$\mathcal F_{\mathrm{cof}}$ of cofinite sets, every member
	of $k\mathcal F$ is infinite. Put
	\[
	RP_{\mathcal F}=\bigcap_{S\in k\mathcal F}RP_S.
	\]
	In particular, $RP_{\mathcal F_{\mathrm{cof}}}=\mathcal C(T)$, where
	$\mathcal C(T)$ is defined in \eqref{eq:localrelation}.
	
	We first compare the sets $C_\varepsilon(U,V)$ with
	finite intersections of return-time sets.
	
	\begin{prop}\label{prop:timecomparison}
		Let $(X,T)$ be minimal and weakly mixing.
		For every finite collection of pairs of nonempty open sets
		$(U_i,V_i)$, $1\le i\le r$, there exist nonempty open
		sets $A,B$ and $\delta>0$ such that
		\begin{equation}\label{eq:timecomparison}
			C_\delta(A,B)\subset\bigcap_{i=1}^rN(U_i,V_i).
		\end{equation}
	\end{prop}

\begin{proof}
If $X$ is a singleton, the conclusion is immediate.
Let $H=\bigcap_{i=1}^rN(U_i,V_i)$, and enumerate
all $C_\varepsilon(P,Q)$, where $P,Q$ belong to a
countable base and $\varepsilon>0$ is rational, as
$C_1,C_2,\ldots$.

We first show that for any $q\ge1$,
\begin{equation}\label{eq:finiteintersection}
C_\delta(A,B)\subseteq
\bigcap_{\ell=1}^q C_{\varepsilon_\ell}(P_\ell,Q_\ell).
\end{equation}
for some nonempty open $A,B$ and $\delta>0$.
By weak mixing, choose a transitive point
$(a,b)\in X^2$ and $k_\ell\ge0$ such that
$(T^{k_\ell}a,T^{k_\ell}b)\in P_\ell\times Q_\ell$.
Choose neighborhoods $A\ni a$, $B\ni b$ with
$T^{k_\ell}A\subset P_\ell$ and
$T^{k_\ell}B\subset Q_\ell$ for all $\ell$.
By uniform continuity, choose $\delta>0$ such that
$\rho(u,v)<\delta$ implies
$\rho(T^{k_\ell}u,T^{k_\ell}v)<\varepsilon_\ell$
for all $\ell$. If $n\in C_\delta(A,B)$, take
$x\in A$, $y\in B$ with $\rho(T^nx,T^ny)<\delta$.
Then $\rho(T^{n+k_\ell}x,T^{n+k_\ell}y)
<\varepsilon_\ell$, proving the inclusion.

Next, for every $M\ge0$, there exists $m$ such that
$C_m\cap\{0,\ldots,M\}=\varnothing$.
Indeed, since $T$ is surjective and $X$ is nontrivial,
choose $p,q\in X$ with $T^Mp\ne T^Mq$. Then
$\eta=\min_{0\le n\le M}\rho(T^np,T^nq)>0$.
By continuity, there are basic neighborhoods $P,Q$
of $p,q$ with $\dist(T^nP,T^nQ)>\eta/2$
for $0\le n\le M$. Thus any rational
$0<\varepsilon<\eta/2$ gives
$C_\varepsilon(P,Q)\cap\{0,\ldots,M\}=\varnothing$.

Suppose that $\bigcap_{\ell=1}^q C_\ell\not\subseteq H$
for every $q\ge1$. We inductively choose
$0=n_0<n_1<\cdots$ such that
\[
n_j\in\bigcap_{\ell=1}^j C_\ell\setminus H.
\]
Indeed, given $n_{j-1}$, choose $m_j$ with
$C_{m_j}\cap\{0,\ldots,n_{j-1}\}=\varnothing$,
and take
$n_j\in\bigcap_{\ell=1}^j C_\ell\cap C_{m_j}
\setminus H$.

We verify condition (2) of Theorem~\ref{thm:criterion}
along $\{n_j:j\ge1\}$. For any nonempty open $U,V$
and $\varepsilon>0$, choose
$C_\ell=C_\eta(P,Q)$ with $P\subset U$,
$Q\subset V$ and $0<\eta<\varepsilon$.
Since $n_j\in C_\ell$ for $j\ge\ell$, we have
$\dist(T^{n_j}U,T^{n_j}V)<\varepsilon$
for all $j\ge\ell$. Thus
\[
\dist(T^{n_j}U,T^{n_j}V)\longrightarrow0
\qquad(j\to\infty).
\]
By the implication $(2)\Rightarrow(1)$ in
Theorem~\ref{thm:criterion}, $(X,T)$ is mixing
along $\{n_j:j\ge1\}$. Hence $n_j\in H$
for all sufficiently large $j$, a contradiction.

Therefore $\bigcap_{\ell=1}^q C_\ell\subseteq H$
for some $q\ge1$. 
By \eqref{eq:finiteintersection}, there exist nonempty
open sets $A,B$ and $\delta>0$ such that
$C_\delta(A,B)\subseteq H$, as required.

\end{proof}

	We also use the following measure-theoretic property
	of closed forward-invariant relations. The construction of an
	invariant pseudometric from their sections goes back
	to McMahon \cite{McMahon78}; see also
	\cite[Theorem~2.3]{GW06}.
	
	\begin{lem}\label{lem:nullsections}
		Let $(X,T)$ be minimal and weakly mixing, and let
		$R\subsetneq X^2$ be closed with
		$(T\times T)R\subset R$. Then
		\begin{equation}\label{eq:nullsections}
			\mu(R[x])=0
			\qquad(x\in X,\ \mu\in M_T(X)).
		\end{equation}
		Moreover,
		\begin{equation}\label{eq:uniformnull}
			\lim_{N\to\infty}\sup_{x,y\in X}
			\frac1N
			\bigl|\{0\le j<N:(x,T^jy)\in R\}\bigr|=0.
		\end{equation}
	\end{lem}
	
\begin{proof}
	Fix $\mu\in M_T(X)$ and put $g(x)=\mu(R[x])$.
	If $x_j\to x$, closedness of $R$ gives
	$\limsup_j\mathbf 1_{R[x_j]}\le\mathbf 1_{R[x]}$ pointwise.
	By the reverse Fatou inequality, $g$ is upper semicontinuous.
	Forward invariance gives $R[x]\subset T^{-1}R[Tx]$, so
	$g(x)\le g(Tx)$ by invariance of $\mu$. The set of maximizers of
	$g$ is nonempty, closed and forward invariant, hence equals $X$.
	Thus all sections have the same measure $c$ and
	$R[x]=T^{-1}R[Tx]\pmod\mu$.

	For $x_j\to x$, the indicators of $R[x_j]\setminus R[x]$
	converge pointwise to zero, again by closedness. Since the sections
	have equal measure, bounded convergence gives
	\[
	\mu(R[x_j]\triangle R[x])
	=2\mu(R[x_j]\setminus R[x])\longrightarrow0.
	\]
	Consequently,
	\[
	d_R(x,y)=\mu(R[x]\triangle R[y])
	\]
	is a continuous pseudometric, invariant under
	$T\times T$ by the preceding equality modulo $\mu$.
	Weak mixing forces $d_R$ to be constant, hence zero.
	Choose a countable dense set $D\subset X$.
	Then $B=\bigcap_{x\in D}R[x]=\bigcap_{x\in X}R[x]$
	is closed and has measure $c$.
	Surjectivity and forward invariance imply $TB\subset B$.
	If $c>0$, minimality gives $B=X$, contradicting
	$R\ne X^2$. Thus $c=0$, proving \eqref{eq:nullsections}.

	Finally, suppose \eqref{eq:uniformnull} fails.
	Then there exist $\eta>0$, $N_k\to\infty$ and
	$x_k,y_k\in X$ such that
	$(\delta_{x_k}\otimes\nu_k)(R)\ge\eta$, where
	$\nu_k=N_k^{-1}\sum_{j=0}^{N_k-1}\delta_{T^jy_k}$.
	Passing to a subsequence, $x_k\to x$ and
$\nu_k\to\nu\in M_T(X)$. Consequently,
$\delta_{x_k}\otimes\nu_k\to\delta_x\otimes\nu$
weakly. Since $R$ is closed, the Portmanteau
inequality gives
\[
\nu(R[x])
=(\delta_x\otimes\nu)(R)
\ge\limsup_k(\delta_{x_k}\otimes\nu_k)(R)
\ge\eta,
\]
contrary to \eqref{eq:nullsections}.
\end{proof}

\subsection{Mixing criteria}

	\begin{thm}\label{thm:family}
		Let $(X,T)$ be minimal, and let $\mathcal F$ be a proper
		Furstenberg family containing $\mathcal F_{\mathrm{cof}}$.
		Then the following conditions are equivalent:
		\begin{enumerate}
			\item $(X,T)$ is $\mathcal F$-mixing.
			\item $C_\varepsilon(U,V)\in\mathcal F$ for all nonempty
			open $U,V\subset X$ and every $\varepsilon>0$.
			\item $RP_{\mathcal F}=X^2$.
			\item For every $S\in k\mathcal F$ and every $x\in X$,
			$P_S(x)$ is a dense $G_\delta$ subset of $X$.
		\end{enumerate}
		If $(X,T)$ is also weakly mixing, these conditions
		are further equivalent to:
		\begin{enumerate}
			\setcounter{enumi}{4}
			\item For every $S\in k\mathcal F$, there exists
			$x_S\in X$ such that
			\begin{equation}\label{eq:onecell}
				\operatorname{int}\overline{P_S(x_S)}
				\ne\varnothing.
			\end{equation}
			\item For every $S\in k\mathcal F$, there exist
			$x_S\in X$ and $\mu_S\in M_T(X)$ such that
			\begin{equation}\label{eq:positivecell}
				\mu_S\bigl(\overline{P_S(x_S)}\bigr)>0.
			\end{equation}
		\end{enumerate}
	\end{thm}
	\begin{rem}

	A family $\mathcal F$ is \emph{full} if $F\cap S$ is infinite
	for every $F\in\mathcal F$ and $S\in k\mathcal F$.
	For full Furstenberg families, Huang, Shao and Ye \cite[Proposition~3.3]{HSY04}
proved that a minimal system is $\mathcal F$-mixing
if and only if $RP_{\mathcal F}^{(n)}=X^n$
for every $n\ge2$. For filters, they also proved the
equivalence of (1) and (4) above \cite[Theorem~3.4]{HSY04}.
Theorem~\ref{thm:family} requires only $RP_{\mathcal F}=X^2$ and
applies to every proper family containing the cofinite sets.
	\end{rem}
	\begin{proof}
		$(1)\Rightarrow(2)$.
		Choose a nonempty open set $W$ with
		$\diam W<\varepsilon$. Then
		\[
		N(U,W)\cap N(V,W)\subset C_\varepsilon(U,V),
		\]
		and the left-hand side belongs to $\mathcal F$.
	
		$(2)\Leftrightarrow(3)$.
By definition, $RP_{\mathcal F}=X^2$ if and only if
$C_\varepsilon(U,V)\cap S\ne\varnothing$ for every
$S\in k\mathcal F$, all nonempty open $U,V$, and
every $\varepsilon>0$. This is equivalent to (2), since
$k(k\mathcal F)=\mathcal F$. Indeed, if
$A\notin\mathcal F$, upward heredity implies
$\mathbb N_0\setminus A\in k\mathcal F$, which
establishes the nontrivial inclusion
$k(k\mathcal F)\subset\mathcal F$.
		
		$(2)\Rightarrow(1)$.
		We first prove weak mixing. Given nonempty open
		$U_1,U_2,V_1,V_2$, choose $k\ge0$ such that
		\[
		W=V_1\cap T^{-k}V_2\ne\varnothing,
		\]
		and a nonempty open $W_0$ with
		$\overline{W_0}\subset W$.
		By minimality, $		X=\bigcup_{r=0}^R T^{-r}W_0$
		for some $R$. Choose $\delta>0$ such that,
		for $0\le r\le R$,
		\[
		\rho(p,q)<\delta,\quad T^rp\in W_0
		\quad\Longrightarrow\quad T^rq\in W.
		\]
		Since $\mathcal F$ is proper, condition~(2) implies
$C_\delta(U_1,T^{-k}U_2)\ne\varnothing$.
Thus there exist $n\ge0$, $a\in U_1$, and
$b\in T^{-k}U_2$ such that
$\rho(T^na,T^nb)<\delta$.
		For some $r\le R$, we have
		$T^{n+r}a\in W_0$ and $T^{n+r}b\in W$.
		Since $  W
=V_1\cap T^{-k}V_2$, we have
\[
T^{n+r}a\in V_1,\qquad
T^{n+r}(T^kb)=T^k(T^{n+r}b)\in V_2.
\]
As $a\in U_1$ and $T^kb\in U_2$, it follows that
\[
n+r\in N(U_1,V_1)\cap N(U_2,V_2).
\]
Thus $(X,T)$ is weakly mixing.
		
		By Proposition~\ref{prop:timecomparison}, there are nonempty open sets
		$A,B$ and $\delta>0$ such that
		\[
		C_\delta(A,B)\subseteq\bigcap_{i=1}^2 N(U_i,V_i).
		\]
		By (2), $C_\delta(A,B)\in\mathcal F$. Since $\mathcal F$ is upward
		hereditary, the intersection also belongs to $\mathcal F$.
		Hence $(X,T)$ is $\mathcal F$-mixing.
		
		$(1)\Rightarrow(4)$.
 Fix $x\in X$ and $S\in k\mathcal F$.
For $m\ge1$, put
\[
G_m=\bigcup_{n\in S}
\{y\in X:\rho(T^nx,T^ny)<1/m\}.
\]
To show that $G_m$ is dense, let $U\subset X$ be
nonempty and choose a finite open cover
$V_1,\ldots,V_r$ with $\diam V_i<1/m$.
By (1), (2), and Proposition~\ref{prop:timecomparison}, $\bigcap_{i=1}^rN(U,V_i)\in\mathcal F$. Thus
there exists
$n\in S\cap\bigcap_{i=1}^rN(U,V_i)$.
Choosing $i$ with $T^nx\in V_i$, we find
$y\in U$ with $T^ny\in V_i$, so $y\in G_m$.

Each $G_m$ is therefore open and dense. To see that
$\bigcap_{m\ge1}G_m\subset P_S(x)$, take $y$ in the intersection
and choose $n_m\in S$ with $\rho(T^{n_m}x,T^{n_m}y)<1/m$.
If $(n_m)$ is unbounded, a subsequence tends to infinity. Otherwise
some $n$ occurs infinitely often, so $T^nx=T^ny$ and the two orbits
agree from time $n$ onward. In either case $y\in P_S(x)$, since $S$
is infinite. The reverse inclusion follows from the definition of
$P_S(x)$, and Baire's theorem now gives~(4).

		$(4)\Rightarrow(3)$.
		Density of $P_S(x)$ for every $x$ gives
		$RP_S=X^2$ for every $S\in k\mathcal F$.
		
		Now suppose that $T$ is weakly mixing.
		The implications $(4)\Rightarrow(5)\Rightarrow(6)$
		follow from density and the full support of invariant
		measures in minimal systems.
		
		$(6)\Rightarrow(3)$.
For each $S\in k\mathcal F$, we have
$\overline{P_S(x_S)}\subset RP_S[x_S]$.
By \eqref{eq:positivecell}, $\mu_S(RP_S[x_S])>0$. Since $RP_S$ is
closed and forward invariant, Lemma~\ref{lem:nullsections}
forces $RP_S=X^2$. Hence $RP_{\mathcal F}=X^2$.
	\end{proof}
	\subsection{Proximal cells along a fixed time set}
	We prove Theorem~\ref{cor:celldichotomy}.

	\begin{proof}
    The implications $(1)\Rightarrow(2)\Rightarrow(3)$
    follow from density and the full support of invariant
    measures in a minimal system.

    $(3)\Rightarrow(1)$.
    Since $\overline{P_S(x)}\subset RP_S[x]$,
    Lemma~\ref{lem:nullsections} forces $RP_S=X^2$.
    Let $\mathcal F_S=\{A\subset\N_0:A\cap S\ne\varnothing\}$.
    As $k\mathcal F_S=\{A\subset\N_0:S\subset A\}$,
    we have $RP_{\mathcal F_S}=RP_S=X^2$.
    The implication $(3)\Rightarrow(4)$ of
    Theorem~\ref{thm:family} now gives (1).

    If these equivalent  conditions  fail, then $RP_S\ne X^2$.
    Lemma~\ref{lem:nullsections} and
    $\overline{P_S(x)}\subset RP_S[x]$ yield the
    asserted nullity and the uniform estimate in
    Remark~\ref{rem:proximalextensions}.
\end{proof}
	
	For minimal systems, mild mixing is equivalent to
	$\mathrm{IP}^*$-mixing \cite{HY04}.
	For minimal weakly mixing systems, Theorem~\ref{thm:family}
	therefore characterizes mild mixing by \eqref{eq:onecell} or
	\eqref{eq:positivecell} for every IP-set $S$.
	If mild mixing fails, some IP-set $S$ has $RP_S\ne X^2$;
	Lemma~\ref{lem:nullsections} then gives \eqref{eq:nullsections}
	and \eqref{eq:uniformnull} for $R=RP_S$.

	\subsection{Completion of the proof of Theorem~\ref{thm:local}}

	In Section~\ref{sec:local}, we established the equivalence
of conditions (1)--(3) in Theorem~\ref{thm:local}.
We now prove the remaining implications
$(3)\Rightarrow(4)\Rightarrow(5)\Rightarrow(3)$,
together with the mildly mixing variant in
Remark~\ref{rem:localvariants}, where conditions
(3)--(5) are understood in their modified form.

	\begin{proof}
		$(4)\Rightarrow(5)$ follows from the full support of invariant measures.
		We prove $(5)\Rightarrow(3)\Rightarrow(4)$ simultaneously for the
		mixing and mildly mixing cases.
		
		$(5)\Rightarrow(3)$.
		The invariant metric on the maximal equicontinuous factor gives
		$RP_S\subset R_{\pi_{\mathrm{eq}}}$ for every infinite $S$.
		Apply (5) with $S=\N$, which is an IP-set. Since
		$\overline{P_{\N}(x)}\subset RP_{\N}[x]$, the pushforward of the
		resulting invariant measure to $X_{\mathrm{eq}}$ has a positive atom.
		Invariance makes all points in its orbit have the same positive mass;
		minimality therefore makes $X_{\mathrm{eq}}$ finite.
		By Lemma~\ref{lem:finitefactor}, its fibers $C_i$ form a clopen cyclic
		decomposition, and each $S_i=T^p|_{C_i}$ is minimal and weakly mixing.
		
		Fix a component $C_j$ and an infinite set $A$, or an IP-set $A$ in
		the mildly mixing case. The set $pA$ has the same property.
		By (5), there are $x\in C_i$ and $\mu\in M_T(X)$ such that
		\[
		\mu\bigl(\overline{P_{pA}^T(x)}\bigr)>0.
		\]
		Separation of distinct components gives
		\[
		P_{pA}^T(x)=P_A^{S_i}(x)\subset C_i.
		\]
		The measure $\mu_i=p\mu|_{C_i}$ is an $S_i$-invariant probability
		measure and gives positive measure to $K=\overline{P_A^{S_i}(x)}$.
		Choose $r\ge0$ with $i+r=j$ modulo $p$, and put
		$f=T^r|_{C_i}:C_i\to C_j$. Uniform continuity gives
		\[
		fK\subset\overline{P_A^{S_j}(fx)}.
		\]
		Thus the $S_j$-invariant measure $\nu=f_*\mu_i$ satisfies
		\[
		\nu\bigl(\overline{P_A^{S_j}(fx)}\bigr)
		\ge\nu(fK)\ge\mu_i(K)>0.
		\]
		Since $A$ and $j$ were arbitrary, Theorem~\ref{thm:family}, with
		the cofinite family or the family of IP$^*$-sets, makes every $S_j$
		mixing or mildly mixing, respectively.
		
		$(3)\Rightarrow(4)$.
		Fix $x\in C_i$. In the mixing case, for every infinite $S$ there is
		$0\le r<p$ such that
		\[
		A=\{k\ge0:pk+r\in S\}
		\]
		is infinite. Theorem~\ref{thm:family} makes $P_A^{S_i}(x)$ dense
		in $C_i$, and uniform continuity of $T^r$ gives
		$P_A^{S_i}(x)\subset P_S^T(x)$.
		
		In the mildly mixing case, every IP-set $S$ contains an IP-set
		$S'\subset S\cap p\N$: group successive terms of a generating sequence
		into finite blocks whose sums are divisible by $p$.
		Then $A=S'/p$ is an IP-set. Mild mixing and Theorem~\ref{thm:family}
		make $P_A^{S_i}(x)$ dense in $C_i$, and again
		$P_A^{S_i}(x)\subset P_S^T(x)$.
		
		Distinct components remain uniformly separated under common iterates,
		so in either case $\overline{P_S^T(x)}=C_i$.
		This proves (4) and the cell-closure statement in
		Remark~\ref{rem:localvariants}.
		Since each $S_i$ is weakly mixing, every equicontinuous factor of
		$(X,T)$ is constant on each $C_i$. Thus the cyclic factor is maximal,
		also in the mildly mixing case. Cell density and
		$RP_S\subset R_{\pi_{\mathrm{eq}}}$ now give
		\begin{equation}\label{eq:cellcomponent}
			RP_S=\bigsqcup_iC_i^2,\qquad
			\overline{P_S(x)}=C_i\quad(x\in C_i)
		\end{equation}
		for every infinite $S$ in the mixing case and every IP-set $S$ in the
		mildly mixing case. 
		
		The decomposition is therefore unique up to cyclic
relabeling. Moreover, weak mixing, total minimality,
or connectedness forces $p=1$.

		Finally, if the corresponding decomposition does not exist,
		the contrapositive of $(5)\Rightarrow(3)$ gives the null-cell assertion
		in Remark~\ref{rem:localvariants}.
	\end{proof}

		\section{Questions}\label{sec:questions}
We conclude with three questions suggested by the preceding results.
Throughout this section, all transformations are homeomorphisms.

\subsection{Asymptotic Pairs of Measures}

Let $\mathcal P(X)$ be the space of Borel probability measures on $X$, equipped with the $1$-Wasserstein metric $W_1$ induced by $\rho$, which metrizes the weak-star topology. Let $T_*\mu=\mu\circ T^{-1}$. For $\mu,\nu\in\mathcal P(X)$, the condition $(\mu,\nu)\in\Asy(T_*)$ means that
$$
\int_X f\circ T^n\,d(\mu-\nu)\longrightarrow0
\qquad(f\in C(X)).
$$
For minimal systems, density of this relation implies mixing. Indeed, for nonempty open sets $U,V\subset X$, choose an asymptotic pair $(\mu,\nu)$ with $\mu(U),\nu(V)>3/4$. Every coupling $\lambda_n$ of $T_*^n\mu$ and $T_*^n\nu$ satisfies
$$
\lambda_n(T^nU\times T^nV)
\ge \mu(U)+\nu(V)-1>\frac12.
$$
Consequently,
$$
\frac12\dist(T^nU,T^nV)
\le W_1(T_*^n\mu,T_*^n\nu)\longrightarrow0.
$$
Theorem~\ref{thm:criterion} now gives mixing. We ask whether the converse holds.

\begin{problem}\label{prob:inducedasymptotic}
	If $(X,T)$ is minimal and mixing, must
	$$
	\overline{\Asy(T_*)}=\mathcal P(X)^2?
	$$
\end{problem}
	The answer is affirmative when $\Asy(T)$ is dense in $X^2$, by
	finite atomic approximation.

Suppose $m$ is an invariant probability measure of full support and $(X,m,T)$
is mixing in the measure-theoretic sense. Then the probability measures
$\mu\ll m$ are weak-star dense in $\mathcal P(X)$ and satisfy
$T_*^n\mu\to m$ \cite[Proposition~1.1(3)]{Polo11}.
Lehrer's model theorem \cite{Lehrer87} gives strictly ergodic, topologically
mixing systems $(X,T)$ for which $(X,m,T)$ fails to be mixing, where $m$ is
the unique invariant probability measure.

\subsection{Relative Asymptotic Pairs}

We next consider a relative version of Corollary~\ref{thm:asy}.
For a factor map $\pi:(X,T)\to(Y,S)$, write
\[
R_\pi=\{(x,x')\in X^2:\pi(x)=\pi(x')\}.
\]

Kamiński, Siemaszko and Szymański
\cite[following Proposition~3.17]{KSS08}
raised the following questions.

\begin{problem}\label{prob:relativewm}
	Let $\pi:(X,T)\to(Y,S)$ be a factor map with $(X,T)$ minimal.
	Must every $K$-extension be relatively weakly mixing?
	More generally, does
	\[
	\overline{\Asy(T)\cap R_\pi}=R_\pi
	\]
	imply that $(R_\pi,T\times T)$ is topologically transitive?
\end{problem}

By \cite[Definition~3.15]{KSS08}, $\pi$ is a $K$-extension if there exists a closed
equivalence relation $K\subseteq R_\pi$ such that
\[
(T\times T)K\subseteq K,\qquad
\bigcap_{n\ge0}(T\times T)^nK=\Delta_X,
\]
and
\[
\overline{\bigcup_{n\ge0}(T\times T)^{-n}K}=R_\pi.
\]
Since the compact sets $(T\times T)^nK$ decrease to $\Delta_X$,
every pair in $K$ is asymptotic. The same holds for every pair
that enters $K$ under forward iteration. These pairs form a dense
subset of $R_\pi$ by the last condition, so
$\overline{\Asy(T)\cap R_\pi}=R_\pi$, as in
\cite[Proposition~3.17]{KSS08}.

The classical Sturmian extension gives a negative answer
to both questions. Let $(X,\sigma)$ be a Sturmian subshift
and $\pi:(X,\sigma)\to(\mathbb T,R_\alpha)$ its natural
factor onto an irrational rotation. By
\cite[Example~3.15]{BGKM02}, $\pi$ is one-to-one
except over a single rotation orbit, where each fiber
consists of exactly two points, which are asymptotic
under both $\sigma$ and $\sigma^{-1}$.

Choose distinct $a,b\in\pi^{-1}(y_0)$ from such a fiber.
Then
\[
R_\pi=\Delta_X\cup
\{(\sigma^k a,\sigma^k b),(\sigma^k b,\sigma^k a):
k\in\mathbb Z\}.
\]
Since
\[
\rho(\sigma^k a,\sigma^k b)\longrightarrow0
\qquad (|k|\to\infty),
\]
and $\sigma$ is aperiodic, every off-diagonal point
of $R_\pi$ is isolated. Thus the nonempty open sets
$\{(a,b)\}$ and $\{(\sigma^{-1}a,\sigma^{-1}b)\}$
satisfy
\[
(\sigma\times\sigma)^n\{(a,b)\}
\cap\{(\sigma^{-1}a,\sigma^{-1}b)\}
=\varnothing\qquad(n\ge0).
\]
Hence $(R_\pi,\sigma\times\sigma)$ is not transitive.

Following Szymański's construction \cite[Propositions~2.1 and~2.2]{Szymanski09}, let
\[
K=\Delta_X\cup
\{(\sigma^k a,\sigma^k b),(\sigma^k b,\sigma^k a):
k\ge0\}.
\]
The orbits of $a$ and $b$ are disjoint because $R_\alpha$ is aperiodic,
and the pair $(a,b)$ is asymptotic. Thus $K$ is a closed equivalence relation. Its $n$th forward image
consists of the diagonal and the displayed pairs with $k\ge n$, so
these images decrease to $\Delta_X$. Every off-diagonal pair in
$R_\pi$ with index $k\in\mathbb Z$ enters $K$ after
$\max\{0,-k\}$ iterates. Hence
$\bigcup_{n\ge0}(\sigma\times\sigma)^{-n}K=R_\pi$, and
$\pi$ is a $K$-extension.
Thus the Sturmian extension provides a counterexample
to both questions in Problem~\ref{prob:relativewm}.

The base of this example is not mixing. This suggests
the following question under a stronger assumption on
the factor.

\begin{conj}\label{conj:relativemixing}
	Let $\pi:(X,T)\to(Y,S)$ be a factor map between minimal systems.
	If $(Y,S)$ is mixing and
	\[
	\overline{\Asy(T)\cap R_\pi}=R_\pi,
	\]
	then $(X,T)$ is mixing.
\end{conj}

We give two observations supporting this conjecture.

Let $\pi_{\mathrm{eq}}:X\to X_{\mathrm{eq}}$ be the maximal
equicontinuous factor map of $X$. Since
$\Asy(T)\subseteq R_{\pi_{\mathrm{eq}}}$, the density assumption gives
$R_\pi\subseteq R_{\pi_{\mathrm{eq}}}$. Thus $\pi_{\mathrm{eq}}$
factors through $\pi$. As $Y$ is mixing, $X_{\mathrm{eq}}$
must be trivial, and hence $X$ is weakly mixing.

The conjecture also holds if $\pi$ has a countable fiber.
Indeed, if $X$ were not mixing, there would be nonempty open
sets $U,V\subset X$ and $n_j\to\infty$ such that
$T^{n_j}U\cap V=\varnothing$.
Passing to a subsequence, let $R$ be the Hausdorff limit of
$\{(T^{n_j}x,x):x\in X\}$.
Then $R$ is closed and $(T\times T)$-invariant, and
$R\cap(V\times U)=\varnothing$, so $R\subsetneq X^2$.

Since $\pi$ is semi-open and $Y$ is mixing, for any nonempty
open sets $W\subset Y$ and $O\subset X$,
\[
S^{n_j}\operatorname{int}\pi(O)\cap W\ne\varnothing
\]
for all sufficiently large $j$. Thus there exists
$z_j\in O$ with $\pi(T^{n_j}z_j)\in W$.
Applying this observation to open sets with closures contained in $W$
and $O$, respectively, and passing to a subsequence gives a point of
$(\pi\times\mathrm{id})(R)$ in $W\times O$. This projection is compact
and hence equals $Y\times X$.
Consequently, for any $y_0\in Y$,
\[
X=\bigcup_{x\in\pi^{-1}(y_0)}R[x].
\]
By Lemma~\ref{lem:nullsections}, each $R[x]$ has measure zero
for every invariant probability, contradicting the
countability of $\pi^{-1}(y_0)$.

\subsection{Higher-Order Mixing}

	\begin{problem}\label{prob:threefoldmixing}
		Let $(X,T)$ be a minimal system with UPE of order two. For every choice of nonempty
		open sets $U_0,U_1,U_2\subset X$, does there exist $N\ge1$ such that
		\[
		U_0\cap T^{-a}U_1\cap T^{-(a+b)}U_2\ne\varnothing
		\qquad(a,b\ge N)?
		\]
	\end{problem}
	
	Theorem~\ref{thm:upe} gives the pairwise
	intersections. The question is whether, for every sufficiently large choice
	of the gaps $a$ and $b$, the three sets have a common point. Dekking and Keane \cite{DK78} constructed a
	minimal mixing substitution system which is not mixing of order
	three, illustrating the distinction between twofold and higher-order
	mixing. 	However, their example has zero entropy and thus does not
	apply to the UPE setting. To the best of our knowledge,
	Problem~\ref{prob:threefoldmixing} remains open.

\end{document}